\documentclass[11pt]{amsart}%
\usepackage{eurosym,amscd,amsmath,amsfonts,amssymb,amsthm,color,graphicx,latexsym,geometry,hyperref}
\usepackage{amsmath}
\usepackage{amsfonts}
\usepackage{amssymb}
\usepackage{graphicx}
\usepackage{tikz}%
\providecommand{\U}[1]{\protect\rule{.1in}{.1in}}
\providecommand{\U}[1]{\protect\rule{.1in}{.1in}}
\providecommand{\U}[1]{\protect\rule{.1in}{.1in}}
\providecommand{\U}[1]{\protect\rule{.1in}{.1in}}
\newtheorem{theorem}{Theorem}[section]

\newtheorem{conjecture}{Conjecture}[section]
\theoremstyle{definition}

\begin{document}
\title[Unimodular multilinear forms with small norms on sequence spaces]{On unimodular multilinear forms with small norms on sequence spaces}
\author[D. Pellegrino]{Daniel Pellegrino}
\address{Departamento de Matem\'{a}tica \\
Universidade Federal da Para\'{\i}ba \\
58.051-900 - Jo\~{a}o Pessoa, Brazil.}
\email{pellegrino@pq.cnpq.br and dmpellegrino@gmail.com}
\author[D. Serrano]{Diana Serrano-Rodr\'{\i}guez}
\address{Departamento de Matem\'{a}ticas\\
\indent Universidad Nacional de Colombia\\
\indent111321 - Bogot\'{a}, Colombia}
\email{dmserrano0@gmail.com and diserranor@unal.edu.co}
\author[J.. Silva]{Janiely Silva}
\address{Departamento de Matem\'{a}tica \\
Universidade Federal da Para\'{\i}ba \\
58.051-900 - Jo\~{a}o Pessoa, Brazil.}
\email{janielymrsilva@gmail.com}
\thanks{D. Pellegrino is supported by CNPq Grant 307327/2017-5 and Grant 2019/0014 Paraiba State Research Foundation (FAPESQ) and J. Silva is supported by CAPES}

\begin{abstract}
The Kahane--Salem--Zygmund inequality is a probabilistic result that
guarantees the existence of special matrices with entries $1$ and $-1$
generating unimodular $m$-linear forms $A_{m,n}:\ell_{p_{1}}^{n}\times
\cdots\times\ell_{p_{m}}^{n}\longrightarrow\mathbb{R}$ (or $\mathbb{C}$) with
relatively small norms. The optimal asymptotic estimates for the smallest
possible norms of $A_{m,n}$ when $\left\{  p_{1},...,p_{m}\right\}
\subset\lbrack2,\infty]$ and when $\left\{  p_{1},...,p_{m}\right\}
\subset\lbrack1,2)$ are well-known and in this paper we obtain the optimal
asymptotic estimates for the remaining case: $\left\{  p_{1},...,p_{m}%
\right\}  $ intercepts both $[2,\infty]$ and $[1,2)$. In particular we prove
that a conjecture posed by Albuquerque and Rezende is false and, using a
special type of matrices that dates back to the works of Toeplitz, we also answer a problem posed by the
same authors. \ 

\end{abstract}
\maketitle

\section{Introduction}

Let $\mathbb{K}$ be the real or complex scalar field. The
Kahane--Salem--Zygmund inequality (see \cite{ab, boas}) asserts that for
positive integers $m,n$ and $p_{1},...,p_{m}\in\lbrack2,\infty]$, there exist
a universal constant $C$ (depending only on $m$), a choice of signs $1$ and
$-1$, and an $m$-linear form $A_{m,n}:\ell_{p_{1}}^{n}\times\cdots\times
\ell_{p_{m}}^{n}\longrightarrow\mathbb{K}$ of the type
\[
A_{m,n}(z^{(1)},...,z^{(m)})=\sum_{j_{1},...,j_{m}=1}^{n}\pm z_{j_{1}}%
^{(1)}\cdots z_{j_{m}}^{(m)},
\]
such that
\[
\Vert A_{m,n}\Vert\leq Cn^{\frac{m+1}{2}-\left(  \frac{1}{p_{1}}+\cdots
+\frac{1}{p_{m}}\right)  }.
\]
An interpolation argument shows that if $p_{1},...,p_{m}\in\lbrack1,2]$, there
is a universal constant $C$ (depending only on $m$), and an $m$-linear form as
above such that
\[
\Vert A_{m,n}\Vert\leq Cn^{1-\frac{1}{\max\{p_{1},...,p_{m}\}}}.
\]
The above estimate appears is essence in Bayart's paper \cite{bayart}. Both
the multilinear and polynomial versions of the Kahane--Salem--Zygmund
inequalities play a fundamental role in modern Analysis (see, for instance,
\cite{ab, bk, vel} and the references therein). However, to the best of the
authors' knowledge, despite the existence of more involved abstract
generalizations of the Kahane--Salem--Zygmund inequality (see \cite{mas}), the
best estimate (i.e., the smallest possible exponent for $n$) for the general
case ($p_{1},...,p_{m}\in\lbrack1,\infty]$) of sequence spaces is still
unknown. Recently, Albuquerque and Rezende (\cite{alb}) have proved that, for
$p_{1},...,p_{m}\in\lbrack1,\infty],$ there is a universal constant $C$
(depending only in $m$) and an $m$-linear form as above satisfying%
\[
\Vert A_{m,n}\Vert\leq Cn^{1-\frac{1}{\gamma}+\sum_{k=1}^{m}\max\left\{
\frac{1}{\gamma}-\frac{1}{p_{k}},0\right\}  },
\]
with
\[
\gamma:=\min\left\{  2,\max\{p_{k}:p_{k}\leq2\}\right\}  .
\]
Note that this last estimate encompasses the previous ones. In this note we
obtain the optimal solution to the general case:

\begin{theorem}
\label{888}Let $m,n$ be positive integers and $p_{1},...,p_{m}\in\left[
1,\infty\right]  $. Then there exist a universal constant $C$ (depending only
on $m$), a choice of signs $1$ and $-1$ and an $m$-linear form $A_{m,n}%
:\ell_{p_{1}}^{n}\times\cdots\times\ell_{p_{m}}^{n}\longrightarrow\mathbb{K}$
of the type
\[
A_{m,n}(z^{(1)},...,z^{(m)})=\sum_{j_{1},...,j_{m}=1}^{n}\pm z_{j_{1}}%
^{(1)}\cdots z_{j_{m}}^{(m)},
\]
such that%
\begin{equation}
\Vert A_{m,n}\Vert\leq Cn^{\frac{1}{\min\left\{  \max\{2,p_{k}^{\ast}\right\}
\}}+\sum_{k=1}^{m}\max\left\{  \frac{1}{2}-\frac{1}{p_{k}},0\right\}  },
\label{818}%
\end{equation}
where $p_{k}^{\ast}$ is the conjugate of $p_{k}.$ Moreover, the exponent
$\frac{1}{\min\left\{  \max\{2,p_{k}^{\ast}\right\}  \}}+\sum_{k=1}^{m}%
\max\left\{  \frac{1}{2}-\frac{1}{p_{k}},0\right\}  $ is optimal.
\end{theorem}

\section{The proof}

We begin by recalling the following estimate obtained by Albuquerque and Rezende:

\begin{theorem}
\label{KSZ_gen}(see \cite{alb}) Let $m,n_{1},\ldots,n_{m}$ be positive
integers and $p_{1},\ldots,p_{m}\in\left[  1,\infty\right]  $. Then there
exist a constant $C$ (depending only on $m$), a choice of signs $1$ and $-1$,
and an $m$-linear form $A:\ell_{p_{1}}^{n_{1}}\times\cdots\times\ell_{p_{m}%
}^{n_{m}}\rightarrow\mathbb{K}$ of the form
\[
A\left(  z^{1},\ldots,z^{m}\right)  =\sum_{j_{1}=1}^{n_{1}}\cdots\sum
_{j_{m}=1}^{n_{m}}\pm_{\mathbf{j}}z_{j_{1}}^{1}\cdots z_{j_{m}}^{m},
\]
such that
\[
\Vert A\Vert\leq C\left(  \sum_{k=1}^{m}n_{k}\right)  ^{1-\frac{1}{\gamma}%
}\prod_{k=1}^{m}n_{k}^{\max\left\{  \frac{1}{\gamma}-\frac{1}{p_{k}%
},0\right\}  },
\]
with $\gamma:=\min\left\{  2,\max\{p_{k}:p_{k}\leq2\}\right\}  .$
\end{theorem}

\subsection{Proof of the inequality (\ref{818})}

We shall prove (\ref{818}) following the more general environment of the above
result. We will show that for positive integers $m,n_{1},\ldots,n_{m}$ and
$p_{1},\ldots,p_{m}\in\left[  1,\infty\right]  $, there is a universal
constant (depending only on $m$), and a $m$-linear form $A:\ell_{p_{1}}%
^{n_{1}}\times\cdots\times\ell_{p_{m}}^{n_{m}}\rightarrow\mathbb{K}$ of the
form
\[
A\left(  z^{1},\ldots,z^{m}\right)  =\sum_{j_{1}=1}^{n_{1}}\cdots\sum
_{j_{m}=1}^{n_{m}}\pm z_{j_{1}}^{1}\cdots z_{j_{m}}^{m},
\]
such that
\begin{equation}
\Vert A\Vert\leq C\left(  \sum_{k=1}^{m}n_{k}\right)  ^{\frac{1}{\rho}}%
\prod_{k=1}^{m}n_{k}^{\max\left\{  \frac{1}{2}-\frac{1}{p_{k}},0\right\}  }.
\label{new}%
\end{equation}
with
\[
\rho:=\min\limits_{k}\left\{  \max\{2,p_{k}^{\ast}\right\}  \}.\text{\ }%
\]
If $p_{k}\geq2$, for all $k=1,\cdots,m$, our estimate coincides with the ones
of Theorem \ref{KSZ_gen}. The same happens when $p_{k}<2$ for all $k=1,...,m.$

Finally, let us suppose (with no loss of generality) that $1\leq d<m$, and
$p_{k}\geq2$, for all $k=1,\cdots,d$ and $p_{k}<2$ for $k=d+1,...,m$. Theorem
\ref{KSZ_gen} guarantees the existence of an $m$-linear form $A:\ell_{p_{1}%
}^{n_{1}}\times\cdots\times\ell_{p_{d}}^{n_{d}}\times\ell_{2}^{n_{d+1}}%
\times\cdots\times\ell_{2}^{n_{m}}\rightarrow\mathbb{K}$ such that
\[
\Vert A\Vert_{\mathcal{L}\left(  \ell_{p_{1}}^{n_{1}}\times\cdots\times
\ell_{p_{d}}^{n_{d}}\times\ell_{2}^{n_{d+1}}\times\cdots\times\ell_{2}^{n_{m}%
};\mathbb{K}\right)  }\leq C\left(  \sum_{k=1}^{m}n_{k}\right)  ^{\frac{1}{2}%
}\prod_{k=1}^{m}n_{k}^{\max\left\{  \frac{1}{2}-\frac{1}{p_{k}},0\right\}  }.
\]
On the other hand, for each $k\notin\{1,\cdots,d\},$ by the monotonicity of
the $\ell_{p}$ norms, the restriction of this form to $\ell_{p_{1}}^{n_{1}%
}\times\cdots\times\ell_{p_{d}}^{n_{d}}\times\ell_{p_{d+1}}^{n_{d+1}}%
\times\cdots\times\ell_{p_{m}}^{n_{m}}$ has norm%

\begin{align}
\Vert A\Vert_{\mathcal{L}\left(  \ell_{p_{1}}^{n_{1}}\times\cdots\times
\ell_{p_{d}}^{n_{d}}\times\ell_{p_{d+1}}^{n_{d+1}}\times\cdots\times
\ell_{p_{m}}^{n_{m}};\mathbb{K}\right)  } &  \leq\Vert A\Vert_{\mathcal{L}%
\left(  \ell_{p_{1}}^{n_{1}}\times\cdots\times\ell_{p_{d}}^{n_{d}}\times
\ell_{2}^{n_{d+1}}\times\cdots\times\ell_{2}^{n_{m}};\mathbb{K}\right)
}\label{2s}\\
&  \leq C\left(  \sum_{k=1}^{m}n_{k}\right)  ^{\frac{1}{2}}\prod_{k=1}%
^{m}n_{k}^{\max\left\{  \frac{1}{2}-\frac{1}{p_{k}},0\right\}  }.\nonumber
\end{align}
Note that in this case%
\[
\rho:=\min\limits_{k}\left\{  \max\{2,p_{k}^{\ast}\right\}  \}=2.
\]
Considering $n_{1}=\cdots=n_{m}=n$ we obtain the proof of (\ref{818}).

\subsection{Proof of the optimality}

The optimality of the case $p_{k}\geq2$ for all $k\in\{1,...,m\}$ is
well-known (it is a consequence of the Hardy--Littlewood inequalities) and the
constant involved does not depend on $p_{1},...,p_{m}$.

More precisely, for all unimodular forms we have%
\[
\Vert A\Vert\geq\frac{1}{\left(  \sqrt{2}\right)  ^{m-1}}n^{\frac{1}%
{2}+\left(  \frac{1}{2}-\frac{1}{p_{1}}\right)  +\cdots+\left(  \frac{1}%
{2}-\frac{1}{p_{m}}\right)  }.
\]
It remains only to prove the optimality of the exponents in the case in which
at least one of the $p_{k}$ is smaller than $2$. We shall split the proof in
three cases:

\begin{itemize}
\item First case: $p_{k}<2$, for all $k=1,\cdots,m$

\item Second case: $p_{k}\geq2$ for only one $k\in\left\{  1,\cdots,m\right\}
$.

\item Third case: the complement of the previous cases.
\end{itemize}

The optimality of the first case seems to be folklore, but for the sake of
completeness we shall provide a proof. In the first case the exponent of $n$
is
\[
\frac{1}{\rho}=\frac{1}{\min\limits_{k}\left\{  \max\{2,p_{k}^{\ast}\right\}
\}}=\frac{p_{j}-1}{p_{j}},
\]
where
\[
p_{j}:=\max\limits_{k}p_{k}.
\]
There is no loss of generality in supposing $j=m$. In the second case (we can
also suppose $k=m$), the exponent of $n$ is also $\frac{p_{m}}{p_{m}-1}.$ For
all $m$-linear forms $A:\ell_{p_{1}}^{n}\times\cdots\times\ell_{p_{m}}%
^{n}\rightarrow\mathbb{K},$ we have%
\[
\sup_{j_{1},...,j_{m-1}}\left(  \sum_{j_{m}=1}^{n}\left\vert A\left(
e_{j_{1}},\dots,e_{j_{m}}\right)  \right\vert ^{\frac{p_{m}}{p_{m}-1}}\right)
^{\frac{p_{m}-1}{p_{m}}}\leq\left\Vert A\right\Vert \sup_{\varphi\in
B_{\left(  \ell_{p_{m}}^{n}\right)  ^{\ast}}}\left(  \sum_{j_{m}=1}%
^{n}\left\vert \varphi\left(  e_{j_{m}}\right)  \right\vert ^{\frac{p_{m}%
}{p_{m}-1}}\right)  ^{\frac{p_{m}-1}{p_{m}}}\leq\left\Vert A\right\Vert .
\]
Thus, for all unimodular $m$-linear forms $A:\ell_{p_{1}}^{n}\times
\cdots\times\ell_{p_{m}}^{n}\rightarrow\mathbb{K},$ we have%
\[
\Vert A\Vert\geq n^{\frac{p_{m}-1}{p_{m}}},
\]
and this guarantees the optimality of the exponent for the first and second cases.

It remains to prove the $m$-linear case when at least two $p_{i}\in\left[
2,\infty\right]  $ and at the same time at least one $p_{i}\in\lbrack1,2)$.

We shall proceed by induction on $m$. The case of bilinear forms is completed
by the previous steps. So, let us suppose that the result is valid for
$\left(  m-1\right)  $-linear forms and let us prove for $m$-linear forms. So,
our induction hypothesis is that for all $p_{i}\in\left[  1,\infty\right]  $
and $i=1,...,m-1$ we have (for all unimodular forms $A:\ell_{p_{1}}^{n}%
\times\cdots\times\ell_{p_{m-1}}^{n}\rightarrow\mathbb{K}$)%
\[
\Vert A\Vert\geq D_{m-1}n^{\frac{1}{\min\left\{  \max\{2,p_{1}^{\ast
}\},...,\max\{2,p_{m-1}^{\ast}\}\right\}  }+\sum_{k=1}^{m-1}\max\left\{
\frac{1}{2}-\frac{1}{p_{k}},0\right\}  }%
\]
and we want to prove that (for all unimodular forms $A:\ell_{p_{1}}^{n}%
\times\cdots\times\ell_{p_{m}}^{n}\rightarrow\mathbb{K}$) we have%
\[
\Vert A\Vert\geq D_{m-1}n^{\frac{1}{\min\left\{  \max\{2,p_{1}^{\ast
}\},...,\max\{2,p_{m}^{\ast}\}\right\}  }+\sum_{k=1}^{m}\max\left\{  \frac
{1}{2}-\frac{1}{p_{k}},0\right\}  }.
\]
Recalling that it just remains to prove the case when at least two $p_{i}%
\in\left[  2,\infty\right]  $ and at the same time at least one $p_{i}%
\in\lbrack1,2)$, we have
\[
\rho=\min\limits_{k}\left\{  \max\{2,p_{k}^{\ast}\}\right\}  =2.
\]
So, we shall prove that for all unimodular $m$-linear forms $A:\ell_{p_{1}%
}^{n}\times\cdots\times\ell_{p_{m}}^{n}\rightarrow\mathbb{K}$ (when at least
two $p_{i}\in\left[  2,\infty\right]  $ and at the same time at least one
$p_{i}\in\lbrack1,2)$) we have
\[
\Vert A\Vert\geq D_{m}n^{\frac{1}{2}+\sum_{k=1}^{m}\max\left\{  \frac{1}%
{2}-\frac{1}{p_{k}},0\right\}  }.
\]
We can suppose that $p_{m}\in\lbrack1,2)$. In this case, for any unimodular
$m$-linear form
\[
A:\ell_{p_{1}}^{n}\times\cdots\times\ell_{p_{m}}^{n}\rightarrow\mathbb{K}%
\]
we have, by the Induction Hypothesis,
\begin{align*}
\Vert A\Vert &  \geq\sup\left\{  \left\vert A\left(  x_{j_{1}}^{\left(
1\right)  },\dots,x_{j_{m-1}}^{\left(  m-1\right)  },\left(  1,0,...0\right)
\right)  \right\vert :\text{ }%
{\textstyle\sum\limits_{j_{k}=1}^{n}}
\left\vert x_{j_{k}}^{(k)}\right\vert ^{p_{k}}\leq1\text{ for all }1\leq k\leq
m-1\right\} \\
&  \geq D_{m-1}n^{\frac{1}{2}+\sum\limits_{k=1}^{m-1}\max\left\{  \frac{1}%
{2}-\frac{1}{p_{k}},0\right\}  }\\
&  =D_{m-1}n^{\frac{1}{2}+\sum\limits_{k=1}^{m}\max\left\{  \frac{1}{2}%
-\frac{1}{p_{k}},0\right\}  }.
\end{align*}

\section{The conjecture of Albuquerque--Rezende is false}

The following conjecture was proposed by Albuquerque and Rezende (see \cite[Conjecture 3.3]{alb}):

\begin{conjecture}
Let $p_{1},\dots,p_{m}\in\lbrack1,\infty]$. There exist $B_{m},C_{m}>0$
(depending only on $m$) such that
\begin{equation}
B_{m}\leq\inf\frac{\Vert A\Vert}{\left(  \sum_{k=1}^{m}n_{k}^{1-\frac
{1}{\gamma}}\right)  \cdot\prod_{k=1}^{m}n_{k}^{\max\left\{  \frac{1}{\gamma
}-\frac{1}{p_{k}},0\right\}  }}\leq C_{m}, \label{113}%
\end{equation}
with $\gamma:=\min\left\{  2,\max\{p_{k}:p_{k}\leq2\}\right\}  $, and the
infimum is calculated over all unimodular $m$-linear forms $A:\ell_{p_{1}%
}^{n_{1}}\times\cdots\times\ell_{p_{m}}^{n_{m}}\rightarrow\mathbb{K}$ and the
exponents involved are sharp.
\end{conjecture}

Note that the estimate (\ref{2s}) shows that the conjecture is false. In fact,
for the sake of illustration, let us choose $m=3$, $p_{1}=3/2$ and
$p_{2}=p_{3}=3.$ By (\ref{2s}) there is a universal constant $C$ such that for
all $n_{1},n_{2},n_{3}$ there exist a unimodular trilinear form $A:\ell
_{p_{1}}^{n_{1}}\times\ell_{p_{2}}^{n_{2}}\times\ell_{p_{3}}^{n_{3}%
}\rightarrow\mathbb{K}$ satisfying%
\[
\Vert A\Vert\leq C\left(  n_{1}+n_{2}+n_{3}\right)  ^{1/2}n_{2}^{1/6}%
n_{3}^{1/6}.
\]
Thus, if (\ref{113}) was valid, we would have
\[
0<\frac{C\left(  n_{1}+n_{2}+n_{3}\right)  ^{1/2}n_{2}^{1/6}n_{3}^{1/6}%
}{\left(  n_{1}^{1/3}+n_{2}^{1/3}+n_{3}^{1/3}\right)  n_{2}^{1/3}n_{3}^{1/3}}%
\]
for all $n_{1},n_{2},n_{3}$, and this is impossible.

We end this paper by answering a problem posed in \cite{alb} for
complex-valued versions of the Kahane--Salem--Zygmund inequality. More
precisely, in \cite[Problem 3.6]{alb} the authors ask about the constants
involved in complex-valued versions of the Kahane--Salem--Zygmund inequality,
i.e., when the coefficients $1$ and $-1$ are replaced by complex numbers with
modulo $1$. We shall show that in the bilinear case the former constant can be replaced by
$1$.

Let $p_{1},p_{2}\geq2$ and $n$ such that $n=\max\{n_{1},n_{2}\}$. Borrowing
ideas that date back to Toeplitz \cite{toe} and Littlewood \cite{li} (see also
\cite[page 609]{bh}), we consider a $n\times n$ matrix $(a_{ij})$ defined by
\[
a_{ij}=e^{2\pi i\frac{ij}{n}}.
\]
Note that
\[
\sum_{t=1}^{n}a_{rt}\overline{a_{st}}=n\delta_{rs}.
\]
Define $A:\ell_{p_{1}}^{n_{1}}\times\ell_{p_{2}}^{n_{2}}\rightarrow\mathbb{C}$
by
\[
A\left(  x^{(1)},x^{(2)}\right)  =\sum_{i_{1},i_{2}=1}^{n}a_{i_{1}i_{2}%
}x_{i_{1}}^{(1)}x_{i_{2}}^{(2)}.
\]
Let $x^{(1)}\in B_{\ell_{p_{1}}^{n_{1}}}$ and $x^{(2)}\in B_{\ell_{p_{2}%
}^{n_{2}}}$, where $B_{\ell_{p_{1}}^{n_{1}}}$ and $B_{\ell_{p_{2}}^{n_{2}}}$
are the closed unit balls of $\ell_{p_{1}}^{n_{1}}$ and $\ell_{p_{2}}^{n_{2}}%
$, respectively. Then, completing with zeros, if necessary, consider
$y^{(1)}=(x_{1}^{(1)},\ldots,x_{{n_{1}}}^{(1)},0\ldots,0)$ and $y^{(1)}%
=(x_{1}^{(2)},\ldots,x_{{n_{2}}}^{(2)},0\ldots,0)$ in $B_{\ell_{p_{1}}^{n}}%
$and $B_{\ell_{p_{2}}^{n}}.$ Using the H\"{o}lder inequality, we have%
\begin{align*}
\left\vert A\left(  x^{(1)},x^{(2)}\right)  \right\vert  &  \leq\sum_{i_{2}%
=1}^{n}\left\vert \sum_{i_{1}=1}^{n}a_{i_{1}i_{2}}y_{i_{1}}^{(1)}\right\vert
\left\vert y_{i_{2}}^{(2)}\right\vert \\
&  \leq\left(  \sum_{i_{2}=1}^{n}|y_{i_{2}}^{(2)}|^{2}\right)  ^{\frac{1}{2}%
}\left(  \sum_{i_{2}=1}^{n}\left\vert \sum_{i_{1}=1}^{n}a_{i_{1}i_{2}}%
y_{i_{1}}^{(1)}\right\vert ^{2}\right)  ^{\frac{1}{2}}\\
&  =\left(  \sum_{i_{2}=1}^{n_{2}}\left\vert x_{i_{2}}^{(2)}\right\vert
^{2}\right)  ^{\frac{1}{2}}\left(  \sum_{i_{2}=1}^{n}\left\vert \sum_{i_{1}%
=1}^{n}a_{i_{1}i_{2}}y_{i_{1}}^{(1)}\right\vert ^{2}\right)  ^{\frac{1}{2}}\\
&  \leq\left(  \sum_{i_{2}=1}^{n_{2}}|1|\right)  ^{\frac{1}{2}-\frac{1}{p_{2}%
}}\left(  \sum_{i_{2}=1}^{n_{2}}\left\vert x_{i_{2}}^{(2)}\right\vert ^{p_{2}%
}\right)  ^{\frac{1}{p_{2}}}\left(  \sum_{i_{2}=1}^{n}\left\vert \sum
_{i_{1}=1}^{n}a_{i_{1}i_{2}}y_{i_{1}}^{(1)}\right\vert ^{2}\right)  ^{\frac
{1}{2}}\\
&  \leq n_{2}^{\frac{1}{2}-\frac{1}{p_{2}}}\left(  \sum_{i_{2}=1}%
^{n}\left\vert \sum_{i_{1}=1}^{n}a_{i_{1}i_{2}}y_{i_{1}}^{(1)}\right\vert
^{2}\right)  ^{\frac{1}{2}}.
\end{align*}
Since%
\[
\left(  \sum_{i_{2}=1}^{n}\left\vert \sum_{i_{1}=1}^{n}a_{i_{1}i_{2}}y_{i_{1}%
}^{(1)}\right\vert ^{2}\right)  ^{\frac{1}{2}}=\left(  \sum_{i_{2}=1}^{n}%
\sum_{\substack{i_{1}=1\\j_{1}=1}}^{n}y_{i_{1}}^{(1)}\overline{y_{j_{1}}%
^{(1)}}a_{i_{1}i_{2}}\overline{a_{j_{1}i_{2}}}\right)  ^{\frac{1}{2}}=\left(
\sum_{\substack{i_{1}=1\\j_{1}=1}}^{n}y_{i_{1}}^{1}\overline{y_{j_{1}}^{(1)}%
}\underbrace{\sum_{i_{2}=1}^{n}a_{i_{1}i_{2}}\overline{a_{j_{1}i_{2}}}%
}_{n\delta_{i_{1}j_{1}}}\right)  ^{\frac{1}{2}},
\]
we have
\begin{align*}
\left\vert A\left(  x^{(1)},x^{(2)}\right)  \right\vert  &  \leq n_{2}%
^{\frac{1}{2}-\frac{1}{p_{2}}}\left(  \sum_{\substack{i_{1}=1\\j_{1}=1}%
}^{n}y_{i_{1}}^{(1)}\overline{y_{j_{1}}^{(1)}}n\delta_{i_{1}j_{1}}\right)
^{\frac{1}{2}}\\
&  \leq n^{\frac{1}{2}}n_{2}^{\frac{1}{2}-\frac{1}{p_{2}}}\left(  \sum
_{i_{1}=1}^{n}\left\vert y_{i_{1}}^{(1)}\right\vert ^{2}\right)  ^{\frac{1}%
{2}}\\
&  =n^{\frac{1}{2}}n_{2}^{\frac{1}{2}-\frac{1}{p_{2}}}\left(  \sum_{i_{1}%
=1}^{n_{1}}\left\vert x_{i_{1}}^{(1)}\right\vert ^{2}\right)  ^{\frac{1}{2}}\\
&  \leq n^{\frac{1}{2}}n_{2}^{\frac{1}{2}-\frac{1}{p_{2}}}\left(  \sum
_{i_{1}=1}^{n_{1}}|1|\right)  ^{\frac{1}{2}-\frac{1}{p_{1}}}\left(
\sum_{i_{1}=1}^{n_{1}}\left\vert x_{i_{1}}^{(1)}\right\vert ^{p_{1}}\right)
^{\frac{1}{p_{1}}}\\
&  \leq n^{\frac{1}{2}}n_{1}^{\frac{1}{2}-\frac{1}{p_{1}}}n_{2}^{\frac{1}%
{2}-\frac{1}{p_{2}}}.
\end{align*}
Thus
\[
\left\Vert A\right\Vert \leq n^{\frac{1}{2}}n_{1}^{\frac{1}{2}-\frac{1}{p_{1}%
}}n_{2}^{\frac{1}{2}-\frac{1}{p_{2}}}\leq\left(  n_{1}^{\frac{1}{2}}%
+n_{2}^{\frac{1}{2}}\right)  n_{1}^{\frac{1}{2}-\frac{1}{p_{1}}}n_{2}%
^{\frac{1}{2}-\frac{1}{p_{2}}}.
\]
In \cite{alb} it is proved that%
\[
\inf\frac{\left\Vert A\right\Vert }{\left(  n_{1}^{\frac{1}{2}}+n_{2}%
^{\frac{1}{2}}\right)  n_{1}^{\frac{1}{2}-\frac{1}{p_{1}}}n_{2}^{\frac{1}%
{2}-\frac{1}{p_{2}}}}\leq8\sqrt{2\ln9}\approx16.8.
\]
For the complex case, our result shows that%
\[
\inf\frac{\left\Vert A\right\Vert }{\left(  n_{1}^{\frac{1}{2}}+n_{2}%
^{\frac{1}{2}}\right)  n_{1}^{\frac{1}{2}-\frac{1}{p_{1}}}n_{2}^{\frac{1}%
{2}-\frac{1}{p_{2}}}}\leq\inf\frac{\left\Vert A\right\Vert }{\left(
\max\{n_{1},n_{2}\}\right)  ^{\frac{1}{2}}n_{1}^{\frac{1}{2}-\frac{1}{p_{1}}%
}n_{2}^{\frac{1}{2}-\frac{1}{p_{2}}}}\leq1.
\]

The constant $1$ that we have just obtained is optimal in a certain sense: if we fix, for instance, $n_{1}=1$, then it is simple to see that the infimum on the right-hand-side is precisely $1$.


\begin{thebibliography}{99}                                                                                               %


\bibitem {alb}N. Albuquerque, L. Rezende, Asymptotic estimates for unimodular
multilinear forms with small norms on sequence spaces, to appear in Bull.
Braz. Math. Soc.

\bibitem {bayart}F. Bayart, Maximum modulus of random polynomials. Q. J. Math.
63 (2012), no. 1, 21--39.

\bibitem {ab}F. Bayart, D. Pellegrino, J.B. Seoane-Sep\'{u}lveda, The Bohr
radius of the n-dimensional polydisk is equivalent to $\sqrt{\left(  \log
n\right)  /n}$. Adv. Math. 264 (2014), 726--746.

\bibitem {boas}H.P. Boas, Majorant series. Several complex variables (Seoul,
1998). J. Korean Math. Soc. 37 (2000), no. 2, 321--337.

\bibitem {bk}H.P. Boas, D. Khavinson, Bohr's power series theorem in several
variables. Proc. Amer. Math. Soc. 125 (1997), no. 10, 2975--2979.

\bibitem {bh}H. F. Bohnenblust and E. Hille, On the absolute convergence of
Dirichlet series, Ann. of Math. \textbf{32} (1931), 600--622.

\bibitem {li}J.E. Littlewood, On bounded bilinear forms in an infinite number
of variables, Q J Math, Volume os-1, Issue 1, (1930), 164--174.

\bibitem {mas}M. Masty\l o, R. Szwedek, Kahane-Salem-Zygmund polynomial
inequalities via Rademacher processes. J. Funct. Anal. 272 (2017), no. 11, 4483--4512.

\bibitem {vel}J. Santos, T. Velanga, On the Bohnenblust-Hille inequality for
multilinear forms. Results Math. 72 (2017), no. 1-2, 239--244.

\bibitem {toe}O. Toeplitz, Uber eine bei den Dirichletschen Reihen auftretende
Aufgabe aus der Theorie der Potenzreihen vonunendlichvielen
Ver\"{a}nderlichen, Nachrichten von der K\"{o}niglichen Gesellschaft der
Wissenschaften zu G\"{o}ottingen, 417--432 (1913).
\end{thebibliography}
\end{document}